\documentclass[a4paper,11pt,notitlepage,draft]{article}
\usepackage{amsmath}

\begin{document}

\title{The Exact General Solution of Painlev\'{e}'s Sixth Equation (PVI) and The Exact General Solution of the Navier Stokes Equations with Applications to Boundary Layer Problems}         
\author{Lance Roman-Miller\thanks{The author expresses his deepest appreciation and thanks to Professor P. Broadbridge, Head of School of Engineering and Mathematical Sciences, La Trobe University, for supervising preparation of this paper. Original priority documents for this paper were confidentially submitted to the Head of School in August 2010. Further, thanks and dedications are noted at the end of this paper.  In accordance with the author's obligations as a staff member, La Trobe University's Manager for Research Development and Research Services has been informed the authors work. Additional priority documents have been lodged with the Australian Federal Government.}\\%
Research Associate, School of Engineering and Mathematical Sciences\\
La Trobe University\\
Bundoora Campus, Melbourne, Victoria, Australia\\
Staff Email: L.Roman-Miller@latrobe.edu.au\\
Personal Email: LanceRomanMiller@hotmail.com\\
Surface Mail:  PO Box R813, Royal Exchange, 1225, Sydney, NSW, Australia}
\date{\today}          
\maketitle

\section{Introduction}       

This paper provides the first known exact general solutions of Painlev\'{e}'s sixth equation (PVI) and the exact general solutions of the Navier Stokes equations and Prandtl's boundary layer equations. In all cases the exact solutions will be given as infinite series expansions.
This paper is the result of work conducted over a 20 year period on the exact solution of nonlinear partial differential equations and nonlinear ordinary differential equations.\\

In around 1991 the author first conducted experimental trial series solution work using Waterloo Maple Software and Reduce on the solution of nonlinear differential equations in association with G.\ Smith of the University of Technology Sydney's School of Mathematics, eventually leading to publication of work on the Van der Pol equation\cite{r31}. Additional, work for the author's doctoral dissertation\cite{r29} at the University of Wollongong under the supervision of P.\ Broadbridge and T.\ Marchant led to a further publication on the exact integration of a reduced Fisher's equation, a reduced Blasius equation and the Lorenz model\cite{r33}.\\

\section{Painlev\'{e}'s Sixth Equation, PVI}

According to Sachdev\cite[page~423]{r8}, the problem undertaken by Gambier and Painlev\'{e} was to determine  those equations of the form $w''=f(w',w,z)$, where f is rational in $w$ and $w'$, whose solutions (as opposed to the equations themselves) are free from movable critical points (branch points and essential singularities). Fifty such equations where shown to have solutions having no movable critical points. Of these 50 equations, six in particular required extremely elaborate proofs\cite[page~424]{r8} to show that their solutions had no movable critical points. The solutions to the six equations became known as the Painlev\'{e}
transcendents.\\

While noting that the Painlev\'{e} transcendents are, in general, not expressible in the form of classical transcendents (solutions of linear ordinary differential equations, exponential, circular and elliptic functions)\cite[ see the definition of classical transcendents at page~330]{r4}, Ince suggested the `possible' existence\cite[page~345]{r4}, of special cases that could be expressible in terms of classical transcendents. As noted by Sachdev\cite[page~436]{r8}, Ince's conjecture was confirmed in the 1960s and early 1970s by breakthrough work of mathematicians in Belorussia\cite[pages~436-437]{r8}, \cite[pages~476-477, 888]{r10}, \cite{r12}-\cite{r20},
which lead to the construction of special solutions in terms of Airy, Bessel, Weber-Hermite,
Whittaker and hypergeometric functions\cite[page~442]{r8}. Additional work in the important area of construction of special solutions was, as noted by Sachdev, also undertaken by Fokas and Ablowitz\cite[page~452]{r8}.\\

Ordinary differential equations whose solutions are of the P type (having no movable critical points) have also been seen to have an important relationship to physical properties \cite[pages~427-429]{r8} of waves, in which waves satisfying partial differential equations used to produce the ordinary differential equations interact in such a way that their solutions emerge, after interaction, substantially unchanged.\\

Considerable advances have been made in relation to construction of the exact solutions of PI-PVI. Davis produced a very detailed account of the solutions (both analytic and Laurent) of PI and PII\cite[chapter~8]{r2}. Davis further extended his results by use of the method of analytic continuation\cite[chapter~9]{r2} and provided methods for determination of the location of poles by equating Laurent and analytic expansions in the neighborhood of a singular point\cite[pages~229-232]{r2}. Davis developed the method of analytic continuation both directly (into the complex plane)\cite[pages~263-266 ]{r2} and
also by way of the method of ``pole vaulting''\cite[page~245]{r2} (in which values of $y$ and $y'$ were found on alternate sides of a pole using a Laurent expansion after the location of a pole has been determined \cite[page~239-232]{r2}.\\

Despite an extensive analysis of PI and PII\cite[chapters~7-9]{r2} using Taylor series and Laurent series expansions and also the generation of special solutions for PII-PVI that take the form of classical transcendents  \cite{r12}-\cite{r20}, no exact general solution of PVI, outside the solutions documented in this paper, appear to have been found. Where a solution $w=\phi(z,C_{1},C_{2})$ of $w''=f(w',w,z)$ is said to be the general solution\cite[page~453]{r8} if for appropriate choices of $C_{1}$ and $C_{2}$, $w$ can be shown to be the solution of an initial value problem in the neighborhood of $z_{0}$.\\

\section{The Exact General Analytic Solution of PVI}    

The form of PVI as given by Murphy\cite[page~183]{r6} has been designated for study. PVI can be written as

\begin{eqnarray}\label{e2}
2x^2y\left(x-1\right)^2\left(y-1\right)\left(y-x\right)y''&=&x^2\left(x-1\right)^2\left[3y^2-2y\left(1+x\right)+x\right]p^2 \nonumber\\
&&-2xy \left( x-1 \right) \left( y-1 \right) \left[ \left( 2x-1 \right)y+x^2 \right] p  \\
&&+2 \Big[ \alpha y^2 \left( y-1 \right)^2 \left( y-x \right )^2+
\beta x \left( y-1 \right)^2 \left( y-x \right)^2 \nonumber \\
&&+\gamma y^2 \left(x-1 \right) \left( y-x \right)^2+
\delta x y^2 \left(x-1 \right) \left( y-1 \right )^2   \Big], \nonumber
\end{eqnarray}

where in (\ref{e2}) $p=y'$ and $q=y''$ and coefficients $a_{0}$, $a_{1}$, $a_{2}$ and $a_{3}$ as appearing in Murphy have been replaced by $\alpha$, $\beta$, $\gamma$ and $\delta$ respectively.\\

In order to avoid the singularities at 0 and 1 in equation (\ref{e2}) the coordinate system will be shifted to $-1$, by way of the transformation $x=x'-1$ {\em after which primes are dropped}, and $x'$ replaced simply by $x$. After the transformation $x=x'-1$ is effected, equation (\ref{e2}) transforms to the following equation:

\begin{equation}\label{e4}
\begin{array}{l}
y^3y''[2x^4-12x^3+26x^2-24x+8]\\
+y^2y''[-2x^5+12x^4-26x^3+24x^2-8x]\\
+yy''[2x^5-14x^4+38x^3-50x^2+32x-8]\\
=y^2p^2[3x^4-18x^3+39x^2-36x+12]\\
+yp^2[-2x^5+12x^4-26x^3+24x^2-8x]\\
+p^2[x^5-7x^4+19x^3-25x^2+16x-4]\\
+y^3p[-4x^3+18x^2-26x+12]\\
+y^2p[-2x^4+14x^3-36x^2+40x-16]\\
+yp[2x^4-10x^3+18x^2-14x+4]\\
+2 \alpha y^6\\
-4 \alpha xy^5\\
+y^4[x^2(2\alpha +2\delta )+x(4\alpha +2\beta +2\gamma -6\delta)+(-4\alpha-2\beta-4\gamma +4\delta )]\\
+y^3[x^2(-4\alpha -4\beta -4\gamma -4\delta)+x(4\alpha +4\beta +12\gamma +12\delta )+(-8\gamma -8\delta )]\\
+y^2[x^3(2\beta +2\gamma )+x^2(2\alpha +2\beta-8\gamma +2\delta )+x(-4\alpha -8\beta +10\gamma -6\delta )\\
\ \ \ +(2\alpha +4\beta -4\gamma +4\delta )]\\
+y[-4\beta x^3+8\beta x^2-4\beta x]\\
+[2\beta x^3-6\beta x^2+6\beta x-2\beta ].
\end{array}
\end{equation}\\

The next body of calculations becomes somewhat involved and space does not permit the setting out of each element of the calculation. However, part of the first member (the term commencing with $y^3y''$) will be considered in detail, to illustrate the general principles involved. \\

Let $y=\sum_{i=0}^{\infty}a_{i}x^i$, then\\

\begin{eqnarray}
y^2&=&\sum_{i=0}^{\infty}a_{i}x^i \sum_{j=0}^{\infty}a_{j}x^j\nonumber\\
&=&\sum_{i=0}^{\infty}\sum_{j=0}^\infty a_{i}a_{j}x^{i+j}.\label{e6}
\end{eqnarray}

The next step in equation (\ref{e6}) is to let $i+j=j'$ then\\

\begin{eqnarray}\label{e8}
y^2&=&\sum_{i=0}^\infty \sum_{j'=i}^\infty a_{i}a_{j'-i}x^{j'}.
\end{eqnarray}

Primes are now dropped in (\ref{e8}) and the order of summation, that is summation by columns is used to replace summation by rows. Accordingly equation (\ref{e8}) now becomes\\

\begin{eqnarray}
y^2&=&\sum_{j=0}^{\infty}\sum_{i=0}^{j}a_{i}a_{j-i}x^j.\label{e10}
\end{eqnarray}

The practical application of the theory of alternation of summation by rows and columns is presented in detail in examples considered in Schwatt\cite{r22} together with additional examples of operational aspects of series including substantial detail on general expressions for higher derivatives of functions and associated Maclaurin series expansions. Additional aspects of the theory of infinite series that provide theoretical support for the operational aspects of the theory of infinite series can be found in Knopp\cite{r24} and in Bromwich\cite{r26}. Comtet\cite{r25} also provides a contemporary account of combinatorics and its applications to infinite series.\\

The procedure for the development of $y^3$ will now be given. It will then be obvious that the method applied to $y^3$ can then be systematically developed for higher order differential coefficients. From the result for $y^2$ in equation (\ref{e10}), $y^3$ can be written as

\begin{eqnarray}
y^3 &= &\sum_{j=0}^{\infty}\sum_{i=0}^{j}a_{i}a_{j-i}x^j \sum_{k=0}^{\infty}a_{k}x^k\nonumber\\
&=&\sum_{j=0}^{\infty} \sum_{i=0}^{j} \sum_{k=0}^{\infty} a_{i} a_{j-i} a_{k} x^{j+k}      \text{\quad let $j+k=k'$ so $k=k'-j$} \nonumber\\
&=&\sum_{j=0}^{\infty} \sum_{i=0}^{j} \sum_{k'=j}^{\infty} a_{i} a_{j-i} a_{k'-j} x^{k'}\nonumber\\
&=&\sum_{j=0}^{\infty} \sum_{i=0}^{j} \sum_{k=j}^{\infty} a_{i} a_{j-i} a_{k-j} x^{k}\text{\quad  on dropping primes} \nonumber\\
&=&\sum_{k=0}^{\infty} \sum_{i=0}^{j} \sum_{j=0}^{k} a_{i} a_{j-i} a_{k-j} x^{k} \text{ \quad  on changing the order of summation} \nonumber\\
{\rm \ \ so, \ \ } y^3&=&\sum_{k=0}^{\infty} \sum_{j=0}^{k} \sum_{i=0}^{j} a_{i} a_{j-i} a_{k-j} x^{k} \label{e12}
\end{eqnarray}

Next the expression $y^3y''$ will be fully developed so as to permit equation of coefficients of integral powers of $x$ i.e.\ of $x^n$,\\

\begin{eqnarray}
y^3y''&=&\sum_{p=0}^{\infty}\sum_{q=0}^{p}\sum_{r=0}^{q}a_{r}a_{q-r}a_{p-q}x^p\sum_{s=0}^{\infty}(s+2)(s+1)a_{s+2}x^s\nonumber\\
&=&\sum_{p=0}^{\infty}\sum_{q=0}^{p}\sum_{r=0}^{q}\sum_{s=0}^{\infty}(s+2)(s+1)a_{r}a_{q-r}a_{p-q}a_{s+2}x^{p+s}\nonumber\\
&=&\sum_{p=0}^{\infty}\sum_{q=0}^{p}\sum_{r=0}^{q}\sum_{s=p}^{\infty}(s-p+2)(s-p+1)a_{r}a_{q-r}a_{p-q}a_{s-p+2}x^{s}\nonumber\\
{\rm \ \ so,\ \ } y^3y'' &=&\sum_{s=0}^{\infty}\sum_{p=0}^{s}\sum_{q=0}^{p}\sum_{r=0}^{q}(s-p+2)(s-p+1)a_{r}a_{q-r}a_{p-q}a_{s-p+2}x^{s}.\label{e14}
\end{eqnarray}

The full term associated with the differential coefficient $y^3y''$ in equation (\ref{e4}) is $y^3y''[2x^4-12x^3+26x^2-24x+8]$, where the expression $2x^4$ in this term can be satisfactorily resolved in association with $y^3y''$ to yield an expression in respect of which coefficients of $x^i$ can simply be equated. Resolution of the expression $y^3y''2x^4$ will be considered as representative of the calculations involved. From equation (\ref{e14}) we have, after relabeling dummy indices in (\ref{e14}) as follows
$s->i$, $p->j$ and $q->k$,

\begin{eqnarray}
y^3y''2x^4&=&\sum_{i=0}^{\infty}\sum_{j=0}^{i}\sum_{k=0}^{j}\sum_{r=0}^{k}(i-j+2)(i-j+1)a_{r}a_{k-r}a_{j-k}a_{i-j+2}x^{i}2x^4\label{e16}\\
&=&2\sum_{i=4}^{\infty}\sum_{j=0}^{i-4}\sum_{k=0}^{j}\sum_{r=0}^{k}(i-2-j)(i-3-j)a_{r}a_{k-r}a_{j-k}a_{i-2-j}x^{i}.\label{e18} 
\end{eqnarray}

Where the transition from (\ref{e16}) to (\ref{e18}) was effected using the substitution $i'=i+4$ and then substituting in $i=i'-4$ and subsequently dropping primes.\\ 

In respect of subsequent equation of coefficients of powers of $x^i$ clearly equation (\ref{e18}) is \emph{only defined} for $i \ge 4$. This is an 
\emph{important} point as some subsidiary equations for coefficients are defined only for a fixed value of i and others are defined over a range of values of i, as is customarily observed in calculating solutions to linear ordinary differential equations. The members found from resolving each member of equation (\ref{e4}) in a manner similar to that undertaken above for $y^3y''2x^4$ so as to permit equation of coefficients of $x^i$ are listed below. Members (1)-(16) below arise from the left hand side of equation (\ref{e4}) and members (17) -(65) arise from the right side of equation (\ref{e4}). Upon equating members below, six defining equations need to be satisfied, five from $i=0$ to $i=4$ inclusive and the general equation for $i \ge 5$.\\

\newpage

\begin{center}
\underline{Members (terms) Arising From Equating Coefficients of $ x^{i} $}
\underline{After Substitution of $y=\sum_{i=0}^{\infty}a_{i}x^i$ into Equation (\ref{e4})}
\end{center}%

\begin{equation}\label{e20}
\begin{array}{lll}
+2\sum_{j=0}^{i-4}\sum_{k=0}^{j}\sum_{r=0}^{k}(i-2-j)(i-3-j)a_{r}a_{k-r}a_{j-k}a_{i-2-j}&i\ge4&[1]\\
\pagebreak
-12\sum_{j=0}^{i-3}\sum_{k=0}^{j}\sum_{r=0}^{k}(i-1-j)(i-2-j)a_{r}a_{k-r}a_{j-k}a_{i-1-j}&i\ge3&[2]\\
+26\sum_{j=0}^{i-2}\sum_{k=0}^{j}\sum_{r=0}^{k}(i-j)(i-1-j)a_{r}a_{k-r}a_{j-k}a_{i-j}&i\ge2&[3]\\
-24\sum_{j=0}^{i-1}\sum_{k=0}^{j}\sum_{r=0}^{k}(i+1-j)(i-j)a_{r}a_{k-r}a_{j-k}a_{i+1-j}&i\ge1&[4]\\
+8\sum_{j=0}^{i}\sum_{k=0}^{j}\sum_{r=0}^{k}(i-j+2)(i-j+1)a_{r}a_{k-r}a_{j-k}a_{i-j+2}&i\ge 0&[5]\\
-2\sum_{j=0}^{i-5}\sum_{k=0}^{j}(i-3-j)(i-4-j)a_{k}a_{j-k}a_{i-3-j}&i\ge 5&[6]\\
+12\sum_{j=0}^{i-4}\sum_{k=0}^{j}(i-2-j)(i-3-j)a_{k}a_{j-k}a_{i-2-j}&i\ge 4&[7]\\
-26\sum_{j=0}^{i-3}\sum_{k=0}^{j}(i-1-j)(i-2-j)a_{k}a_{j-k}a_{i-1-j}&i\ge 3&[8]\\
+24\sum_{j=0}^{i-2}\sum_{k=0}^{j}(i-j)(i-1-j)a_{k}a_{j-k}a_{i-j}&i\ge 2&[9]\\
-8\sum_{j=0}^{i-1}\sum_{k=0}^{j}(i+1-j)(i-j)a_{k}a_{j-k}a_{i+1-j}&i\ge 1&[10]\\
+2\sum_{j=0}^{i-5}(i-3-j)(i-4-j)a_{j}a_{i-3-j}&i\ge 5&[11]\\
-14\sum_{j=0}^{i-4}(i-2-j)(i-3-j)a_{j}a_{i-2-j}&i\ge 4&[12]\\
+38\sum_{j=0}^{i-3}(i-1-j)(i-2-j)a_{j}a_{i-1-j}&i\ge 3&[13]\\
-50\sum_{j=0}^{i-2}(i-j)(i-1-j)a_{j}a_{i-j}&i\ge 2&[14]\\
+32\sum_{j=0}^{i-1}(i+1-j)(i-j)a_{j}a_{i+1-j}&i\ge 1&[15]\\
-8\sum_{j=0}^{i}(i+2-j)(i+1-j)a_{j}a_{i+2-j}&i\ge 0&[16]\\
+3\sum_{j=0}^{i-4}\sum_{k=0}^{j}\sum_{r=0}^{i-4-j}(r+1)(i-3-j-r)a_{k}a_{j-k}a_{r+1}a_{i-3-j-r}&i\ge 4&[17]\\
-18\sum_{j=0}^{i-3}\sum_{k=0}^{j}\sum_{r=0}^{i-3-j}(r+1)(i-2-j-r)a_{k}a_{j-k}a_{r+1}a_{i-2-j-r}&i\ge 3&[18]\\
+39\sum_{j=0}^{i-2}\sum_{k=0}^{j}\sum_{r=0}^{i-2-j}(r+1)(i-1-j-r)a_{k}a_{j-k}a_{r+1}a_{i-1-j-r}&i\ge 2&[19]\\
-36\sum_{j=0}^{i-1}\sum_{k=0}^{j}\sum_{r=0}^{i-1-j}(r+1)(i-j-r)a_{k}a_{j-k}a_{r+1}a_{i-j-r}&i\ge 1&[20]\\
+12\sum_{j=0}^{i}\sum_{k=0}^{j}\sum_{r=0}^{i-j}(r+1)(i+1-j-r)a_{k}a_{j-k}a_{r+1}a_{i+1-j-r}&i\ge 0&[21]\\
-2\sum_{j=0}^{i-5}\sum_{k=0}^{j}(k+1)(j-k+1)a_{k+1}a_{j-k+1}a_{i-5-j}&i\ge 5&[22]\\
+12\sum_{j=0}^{i-4}\sum_{k=0}^{j}(k+1)(j-k+1)a_{k+1}a_{j-k+1}a_{i-4-j}&i\ge 4&[23]\\
-26\sum_{j=0}^{i-3}\sum_{k=0}^{j}(k+1)(j-k+1)a_{k+1}a_{j-k+1}a_{i-3-j}&i\ge 3&[24]\\
+24\sum_{j=0}^{i-2}\sum_{k=0}^{j}(k+1)(j-k+1)a_{k+1}a_{j-k+1}a_{i-2-j}&i\ge 2&[25]\\
-8\sum_{j=0}^{i-1}\sum_{k=0}^{j}(k+1)(j-k+1)a_{k+1}a_{j-k+1}a_{i-1-j}&i\ge 1&[26]\\
+\sum_{j=0}^{i-5}(j+1)(i-4-j)a_{j+1}a_{i-4-j}&i\ge 5&[27]\\
-7\sum_{j=0}^{i-4}(j+1)(i-3-j)a_{j+1}a_{i-3-j}&i\ge 4&[28]\\
+19\sum_{j=0}^{i-3}(j+1)(i-2-j)a_{j+1}a_{i-2-j}&i\ge 3&[29]\\
-25\sum_{j=0}^{i-2}(j+1)(i-1-j)a_{j+1}a_{i-1-j}&i\ge 2&[30]\\
+16\sum_{j=0}^{i-1}(j+1)(i-j)a_{j+1}a_{i-j}&i\ge 1&[31]\\
-4\sum_{j=0}^{i}(j+1)(i-j+1)a_{j+1}a_{i-j+1}&i\ge 0&[32]\\
-4\sum_{j=0}^{i-3}\sum_{k=0}^{j}\sum_{r=0}^{k}(i-2-j)a_{r}a_{k-r}a_{j-k}a_{i-2-j}&i\ge 3&[33]\\
+18\sum_{j=0}^{i-2}\sum_{k=0}^{j}\sum_{r=0}^{k}(i-1-j)a_{r}a_{k-r}a_{j-k}a_{i-1-j}&i\ge 2&[34]\\
-26\sum_{j=0}^{i-1}\sum_{k=0}^{j}\sum_{r=0}^{k}(i-j)a_{r}a_{k-r}a_{j-k}a_{i-j}&i\ge 1&[35]\\
+12\sum_{j=0}^{i}\sum_{k=0}^{j}\sum_{r=0}^{k}(i-j+1)a_{r}a_{k-r}a_{j-k}a_{i-j+1}&i\ge 0&[36]\\
-2\sum_{j=0}^{i-4}\sum_{k=0}^{j}(i-3-j)a_{k}a_{j-k}a_{i-3-j}&i\ge 4&[37]\\
+14\sum_{j=0}^{i-3}\sum_{k=0}^{j}(i-2-j)a_{k}a_{j-k}a_{i-2-j}&i\ge 3&[38]\\
-36\sum_{j=0}^{i-2}\sum_{k=0}^{j}(i-1-j)a_{k}a_{j-k}a_{i-1-j}&i\ge 2&[39]\\
+40\sum_{j=0}^{i-1}\sum_{k=0}^{j}(i-j)a_{k}a_{j-k}a_{i-j}&i\ge 1&[40]\\
-16\sum_{j=0}^{i}\sum_{k=0}^{j}(i-j+1)a_{k}a_{j-k}a_{i-j+1}&i\ge 0&[41]\\
\end{array}
\end{equation}
\begin{equation*}
\begin{array}{lll}
+2\sum_{j=0}^{i-4}(i-3-j)a_{j}a_{i-3-j}&i\ge 4&[42]\\
-10\sum_{j=0}^{i-3}(i-2-j)a_{j}a_{i-2-j}&i\ge 3&[43]\\
+18\sum_{j=0}^{i-2}(i-1-j)a_{j}a_{i-1-j}&i\ge 2&[44]\\
-14\sum_{j=0}^{i-1}(i-j)a_{j}a_{i-j}&i\ge 1&[45]\\
+4\sum_{j=0}^{i}(i+1-j)a_{j}a_{i+1-j}&i\ge 0&[46]\\
+2\alpha \sum_{j=0}^{i}\sum_{k=0}^{j}\sum_{r=0}^{k}\sum_{s=0}^{r}\sum_{t=0}^{s}a_{t}a_{s-t}a_{r-s}a_{k-r}a_{j-k}a_{i-j}&i\ge 0&[47]\\
-4\alpha \sum_{j=0}^{i-1}\sum_{k=0}^{j}\sum_{r=0}^{k}\sum_{s=0}^{r}a_{s}a_{r-s}a_{k-r}a_{j-k}a_{i-1-j}&i\ge 1&[48]\\
+(2\alpha +2\delta )\sum_{j=0}^{i-2}\sum_{k=0}^{j}\sum_{r=0}^{k}a_{r}a_{k-r}a_{j-k}a_{i-2-j}&i\ge 2&[49]\\
+(4\alpha +2\beta+2\gamma -6\delta )\sum_{j=0}^{i-1}\sum_{k=0}^{j}\sum_{r=0}^{k}a_{r}a_{k-r}a_{j-k}a_{i-1-j}&i\ge 1&[50]\\
+(-4\alpha -2\beta -4\gamma +4\delta )\sum_{j=0}^{i}\sum_{k=0}^{j}\sum_{r=0}^{k}a_{r}a_{k-r}a_{j-k}a_{i-j}&i\ge 0&[51]\\
+(-4\alpha -4\beta -4\gamma -4\delta )\sum_{j=0}^{i-2} \sum_{k=0}^{j}a_{k}a_{j-k}a_{i-2-j}&i\ge 2&[52]\\
+(4\alpha +4\beta +12\gamma +12\delta )\sum_{j=0}^{i-1} \sum_{k=0}^{j}a_{k}a_{j-k}a_{i-1-j}&i\ge 1&[53]\\
+( -8\gamma -8\delta )\sum_{j=0}^{i} \sum_{k=0}^{j}a_{k}a_{j-k}a_{i-j}&i\ge 0&[54]\\
+(2\beta+2\gamma)\sum_{j=0}^{i-3}a_{j}a_{i-3-j}&i\ge 3&[55]\\
+(2\alpha +2\beta -8\gamma +2\delta)\sum_{j=0}^{i-2}a_{j}a_{i-2-j}&i\ge 2&[56]\\
+(-4\alpha -8\beta +10\gamma -6\delta )\sum_{j=0}^{i-1}a_{j}a_{i-1-j}&i\ge 1&[57]\\
+(2\alpha +4\beta -4\gamma +4\delta )\sum_{j=0}^{i}a_{j}a_{i-j}&i\ge 0&[58]\\
-4\beta a_{i-3}&i\ge 3&[59]\\
+8\beta a_{i-2}&i\ge 2&[60]\\
-4\beta a_{i-1}&i\ge 1&[61]\\
+2\beta &i= 3&[62]\\
-6\beta &i= 2&[63]\\
+6\beta &i= 1&[64]\\
-2\beta &i= 0&[65]\\
\end{array}
\end{equation*}\\

Next a number of simple observations can be made about members 1-65, hereinafter referred to as $M_{i}$ for $i=1 \dots 65$. It is to be noted that the highest possible subscripts can only be obtained from members 5 and 16 by setting $j=0$. This fact permits the separation of member 5, $M_{5}$, and member 16, $M_{16}$, in set (\ref{e20}) above, into two separate parts so as to put the defining equation for the coefficients arising from substitution of $y=\sum_{i=0}^{\infty} a_{i}x^i$ into (\ref{e4})  into a form that can be used not merely to define the structure of the solution to equation (\ref{e4}) but also to enable the solution to be \emph{constructed iteratively}. This is achieved in detail as follows.\\

Members 5, $(M_{5})$ can be expanded as follows

\begin{equation}\label{e22}
\begin{array}{ll}
8\sum_{j=0}^{i} \sum_{k=0}^{j} \sum_{r=0}^{k}a_{r}a_{k-r}a_{j-k}a_{i-j+2}(i-j+2)(i-j+1) \ \ \ \\
=\\
8a_{0}^3a_{i+2}(i+2)(i+1)+8\sum_{j=1}^{i}\sum_{k=0}^{j}\sum_{r=0}^{k}a_{r}a_{k-r}a_{j-k}a_{i-j+2}(i-j+2)(i-j+1)\\
(i \ge 0) \ \ \ \ \ \ \ \ \ \ \ \ \ \ \ \ \ \ \ \ \ \ \ \ \ (i\ge 1)
\end{array}
\end{equation}

Similarly, member 16, $M_{16}$, appearing in set (\ref{e20}) can be expanded as follows\\

\begin{equation}\label{e24}
\begin{array}{l}
-8\sum_{j=0}^{i}a_{j}a_{i+2-j}(i+2-j)(i+1-j)\\
=\\
-8a_{0}a_{i+2}(i+2)(i+1)-8\sum_{j=1}^{i}a_{j}a_{i+2-j}(i+2-j)(i+1-j)\\
(i \ge 0)\ \ \ \ \ \ \ \ \ \ \ \ \ \ \ \ \ \  \ \ \ \ \ \ \ \ \ (i\ge 1)
\end{array}
\end{equation}\\

Members of set (\ref{e20}) can now be summed and solved for $a_{i+2}$ after member 5, $M_{5}$, and member 16, $M_{16}$, of equation (\ref{e20}) have been expanded in accordance with equations (\ref{e22}) and (\ref{e24}) to yield the following defining equation for the coefficients of the solution to equation (\ref{e4})

\begin{eqnarray}
a_{i+2} &= & {-\sum_{i=1}^{16}M^{*}_i+\sum_{i=17}^{65}M^{*}_{i} \over 8a_{0}(i+2)(i+1)(a_{0}^2-1)} \ \ \  {\rm for \ } \ \ i\ge 0;  \ \ \ \ a_{0} \ne 1,0 \label{e26}
\end{eqnarray}

It is to be noted that the members appearing in equation (\ref{e26}) now have a star over them to emphasize that members 5 and 16 ($M_{15}$), ($M_{16}$) have been modified to the extent that summation now only occurs from $j=1$ onwards. \\

To demonstrate the validity of equation (\ref{e26}) parameters and values of $a_{0}$ and $a_{1}$ have been set as follows $\alpha = \beta = \gamma = \delta = 1$, $a_{0}=2$ and $a_{1}=1$. After, evaluation using equation (\ref{e26}) the following values were found.\\

$a_{0}=2$, $a_{1}=1$, $a_{2}=5/48$, $a_{3}=311/864$, $a_{4}=18725/20736$, $a_{5}=48313/34560$, $a_{6}=17430769/8957952$, $a_{7}=3838061/1451520$, $a_{8}=65037559477/18059231232$,\\ $a_{9}=95777442903929/19503969730560$.\\

These values \emph{have been confirmed} by using Waterloo Maple Software's dsolve routine.\\

Laurent expansions in the neighborhood of a singularity have also been conducted by the author. However, as the the construction of these expansions involve the use of substantially no new principles beyond those recited above no further elaboration on this point will be undertaken.\\

The next step will be to show how the methods used above to solve ordinary differential equations exactly can be generalized to solve partial differential equations exactly. The general system of the Navier Stokes equations will be solved exactly and the methodology derived above will then be applied to solve Prandt's boundary layer equations.\\

\section{The Navier Stokes Equations and Boundary Layer Flow}

Ludwig Prandtl said of the Navier Stokes equations, ``{\em The general differential equations of viscous fluids lead to mathematical difficulties which may be unconquerable for a long time to come.'' ... ``.}\cite[page~104]{r28} \\

This paper presents what is apparently the first set of \emph{general and exact} solutions of the Navier Stokes equations as well as presenting what appears to be the first set of general and exact solutions to the associated boundary layer equations of Prandtl. \\

At large Reynold's number the flow around a streamlined body, 
with \emph {no} flow separation, can be divided into a boundary layer flow confined to a region of thickness $\delta$ around the body (in which viscous effects are important) and a flow field far from the body which closely resembles non viscous flow\cite[page~160]{r30} (where the Reynolds number is given by $Re=\rho U L/ \mu=UL / \nu$ for characteristic velocity U, length L and density $\rho$ and where $\delta \ll L$, coefficient of viscosity $\mu$, and kinematic viscosity $\nu$).\\

Under extreme conditions the velocity profile around the body can become extremely distorted until separation and associated eddying and turbulence occurs\cite[pages~160 and 168]{r30}. In principle separation is characterized\cite[page~113]{r32} by $\left({\partial u / \partial y}\right )_{y=0}=0$. However, in practise separation is far more complex and turbulence can even occur well \emph{before} separation, leading paradoxically to the flow adhering\cite[page~161]{r30} to the body wall even in the presence of an \emph{adverse pressure gradient}.\\

Once the boundary layer equations have been solved ``some'' estimate of the location of the point of separation can be determined. However, it is to be noted that in practise, as recited above, despite 'theoretical' calculations as to the point of separation, other factors come into play, including the fact that the presence of any turbulence in the flow could actually \emph{delay} flow separation. Accordingly, results from theoretical calculations must used in association with experimental results to select appropriate modifications to body structures (typically aircraft wing and ship hull structures) that minimize drag and in the case of an aircraft wing, increase lift.\\

Typically, the modifications involved in boundary layer flow stabilization include the introduction of slats, slots and slits so as to delay the onset of flow separation\cite[page~268]{r32}. A detailed account of which is given in Schlichting under the heading ``Boundary Layer Control''\cite[chapter~XIII]{r32}. By way of comment on flow stabilization techniques, it is to be noted that even \emph{very slight adjustments} to body structure by way of the introduction of slits, slats and other structural changes designed to maintain stable laminar flow, can yield extremely significant increases\cite[page~275]{r32} in lift and also significantly reduce drag. \\

The Navier Stokes equations and the associated boundary layer equations will now be solved.\\

\section{The Exact General Solution of the Navier Stokes Equations}

The Navier Stokes Equations are written as follows (on the assumption that there is no free surface and that pressure is taken relative to static pressure were no flow to exist). The absence of the body force term also implicitly assumes that no body force other than gravitational forces exist, an assumption that is reasonable for most flow structures likely to be encountered outside highly specialized flows that arise in areas including magneto-hydrodynamics.\\

\begin{eqnarray}
u{\partial u \over \partial x}+v{\partial u \over \partial y}+w{\partial u \over \partial z}+{\partial u \over \partial t}&=&-{1 \over \rho}{\partial P \over \partial x}+\nu\nabla^2 u\nonumber\\
u{\partial v \over \partial x}+v{\partial v \over \partial y}+w{\partial v \over \partial z}+{\partial v \over \partial t}&=&-{1 \over \rho}{\partial P \over \partial y}+\nu\nabla^2 v\label{e28}\\
u{\partial w \over \partial x}+v{\partial w \over \partial y}+w{\partial w \over \partial z}+{\partial w \over \partial t}&=&-{1 \over \rho}{\partial P \over \partial z}+\nu\nabla^2 w\nonumber
\end{eqnarray}

Equation (\ref{e28}) is to be solved in conjunction with the continuity equation so as to provide four equations in four unknowns\cite[page~144]{r34} $u$, $v$, $w$ and $P$.\\ 

The continuity equation is given by\\

\begin{eqnarray}
0&=&{\partial u \over \partial x}+{\partial v \over \partial y}+{\partial w \over \partial z}.
\end{eqnarray}

No new principles of solution methodology are introduced beyond those already considered in the solution of PVI. However, for the sake of completeness, the general principle in relation to resolution of nonlinearity for purposes of equating coefficients will again be illustrated in some detail.\\

For $u=\sum_{i=0}^{\infty}\sum_{j=0}^{\infty}\sum_{k=0}^{\infty}\sum_{l=0}^{\infty}A_{i,j,k,l}x^iy^jz^kt^l$ the expression $u{\partial u / \partial x}$ will be resolved so as to permit equating of coefficients.\\

\begin{eqnarray}
{\partial u \over \partial x}&=&{\partial \over \partial x}\sum_{i=0}^{\infty}\sum_{j=0}^{\infty}\sum_{k=0}^{\infty}\sum_{l=0}^{\infty}A_{i,j,k,l}x^iy^jz^kt^l\nonumber\\
&=&\sum_{i=0}^{\infty}\sum_{j=0}^{\infty}\sum_{k=0}^{\infty}\sum_{l=0}^{\infty}iA_{i,j,k,l}x^{i-1}y^jz^kt^l {\rm \ \  let \ \ } i=i'+1  \nonumber\\
{\partial u \over \partial x}&=&\sum_{i=0}^{\infty}\sum_{j=0}^{\infty}\sum_{k=0}^{\infty}\sum_{l=0}^{\infty}(i+1)A_{i+1,j,k,l}x^{i}y^jz^kt^l \qquad \text{after dropping primes}\label{e31}
\end{eqnarray}

In calculating terms including $\partial u /\partial x$ above, terms including $i=-1$ that result in a value of zero associated with $i+1$, that arose during use of the transformation $i=i'+1$, see \eqref{e31}, have obviously not been included. Using \eqref{e31} we obtain,

\begin{eqnarray}
u{\partial u \over \partial x}&=&\sum_{p=0}^{\infty}\sum_{q=0}^{\infty}\sum_{r=0}^{\infty}\sum_{s=0}^{\infty}A_{p,q,r,s}x^py^qz^rt^s \sum_{i=0}^{\infty}\sum_{j=0}^{\infty}\sum_{k=0}^{\infty}\sum_{l=0}^{\infty}(i+1)A_{i+1,j,k,l}x^{i}y^jz^kt^l\nonumber\\
&=&\sum_{p=0}^{\infty}\sum_{q=0}^{\infty}\sum_{r=0}^{\infty}\sum_{s=0}^{\infty}\sum_{i=0}^{\infty}\sum_{j=0}^{\infty}\sum_{k=0}^{\infty}\sum_{l=0}^{\infty}(i+1)A_{p,q,r,s}A_{i+1,j,k,l}x^{i+p}y^{j+q}z^{k+r}t^{l+s}\nonumber\\
&=&\sum_{p=0}^{\infty}\sum_{q=0}^{\infty}\sum_{r=0}^{\infty}\sum_{s=0}^{\infty}\sum_{i=p}^{\infty}\sum_{j=q}^{\infty}\sum_{k=r}^{\infty}\sum_{l=s}^{\infty}(i-p+1)A_{p,q,r,s}A_{i-p+1,j-q,k-r,l-s}x^{i}y^{j}z^{k}t^{l}\nonumber\\
u{\partial u \over \partial x}&=&\sum_{i=0}^{\infty}\sum_{j=0}^{\infty}\sum_{k=0}^{\infty}\sum_{l=0}^{\infty}\sum_{p=0}^{i}\sum_{q=0}^{j}\sum_{r=0}^{k}\sum_{s=0}^{l}(i-p+1)A_{p,q,r,s}A_{i-p+1,j-q,k-r,l-s}x^{i}y^{j}z^{k}t^{l}\label{e30}
\end{eqnarray}

Where as before in the construction of equation (\ref{e30}) the transformations $i=i'-p$, $j=j'-q$, $k=k'-r$ and $l=l'-s$ and the method of changing the direction of summation from rows to columns were used.\\

No new principles other than those recited above have been used in the resolution of other elements of the Navier Stokes equations. Accordingly, setting

\[v=\sum_{i=0}^{\infty}\sum_{j=0}^{\infty}\sum_{k=0}^{\infty}\sum_{l=0}^{\infty}B_{i,j,k,l}x^iy^jz^kt^l\]

\[w=\sum_{i=0}^{\infty}\sum_{j=0}^{\infty}\sum_{k=0}^{\infty}\sum_{l=0}^{\infty}C_{i,j,k,l}x^iy^jz^kt^l\]

and 

\[P=\sum_{i=0}^{\infty}\sum_{j=0}^{\infty}\sum_{k=0}^{\infty}\sum_{l=0}^{\infty}D_{i,j,k,l}x^iy^jz^kt^l\]\\

then the coefficients of the analytic solution of the Navier Stokes equations (\ref{e28}) are defined by the following equations (after resolving each nonlinear member as occurred in relation to $u{\partial u / \partial x}$ so as to permit equating of coefficients of $x^iy^jz^kt^l$)

\begin{equation}\label{e32}
\begin{array}{l}
\sum_{p=0}^{i}\sum_{q=0}^{j}\sum_{r=0}^{k}\sum_{s=0}^{l}\Big((i-p+1)A_{p,q,r,s}A_{i-p+1,j-q,k-r,l-s}\\
+(j-q+1)B_{p,q,r,s}A_{i-p,j-q+1,k-r,l-s}\\
+(k-r+1)C_{p,q,r,s}A_{i-p,j-q,k-r+1,l-s}\Big) +(l+1)A_{i,j,k,l+1}\\
=\\
{-1\over \rho}(i+1)D_{i+1,j,k,l}+\nu\Big((i+2)(i+1)A_{i+2,j,k,l}+(j+2)(j+1)A_{i,j+2,k,l}\\+(k+2)(k+1)A_{i,j,k+2,l}\Big)
\end{array}
\end{equation}\\

\begin{equation}\label{e34}
\begin{array}{l}
\sum_{p=0}^{i}\sum_{q=0}^{j}\sum_{r=0}^{k}\sum_{s=0}^{l}\Big((i-p+1)A_{p,q,r,s}B_{i-p+1,j-q,k-r,l-s}\\
+(j-q+1)B_{p,q,r,s}B_{i-p,j-q+1,k-r,l-s}\\
+(k-r+1)C_{p,q,r,s}B_{i-p,j-q,k-r+1,l-s}\Big) +(l+1)B_{i,j,k,l+1}\\
=\\
{-1\over \rho}(j+1)D_{i,j+1,k,l}+\nu\Big((i+2)(i+1)B_{i+2,j,k,l}+(j+2)(j+1)B_{i,j+2,k,l}\\+(k+2)(k+1)B_{i,j,k+2,l}\Big)
\end{array}
\end{equation}\\

\begin{equation}\label{e36}
\begin{array}{l}
\sum_{p=0}^{i}\sum_{q=0}^{j}\sum_{r=0}^{k}\sum_{s=0}^{l}\Big((i-p+1)A_{p,q,r,s}C_{i-p+1,j-q,k-r,l-s}\\
+(j-q+1)B_{p,q,r,s}C_{i-p,j-q+1,k-r,l-s}\\
+(k-r+1)C_{p,q,r,s}C_{i-p,j-q,k-r+1,l-s}\Big) +(l+1)C_{i,j,k,l+1}\\
=\\
{-1\over \rho}(k+1)D_{i,j,k+1,l}+\nu\Big((i+2)(i+1)C_{i+2,j,k,l}+(j+2)(j+1)C_{i,j+2,k,l}\\+(k+2)(k+1)C_{i,j,k+2,l}\Big).
\end{array}
\end{equation}

The defining equation for the coefficients arising from the continuity equation is

\begin{eqnarray}\label{e38}
0&=&(i+1)A_{i+1,j,k,l}+(j+1)B_{i,j+1,k,l}+(k+1)C_{i,j,k+1,l}
\end{eqnarray}

\section{The Exact General Solution of the Prandtl's\\
 Boundary Layer Equations}

A highly simplified form of the general Navier Stokes equations (\ref{e28}) , known as Prandtl's boundary layer equations (including the continuity equation) is given as follows  \cite[page~109]{r32}

\begin{eqnarray}
u{\partial u \over \partial x}+v{\partial u \over \partial y}+{\partial u \over \partial t}&=&-{1 \over \rho}{\partial P \over \partial x}+\nu\ {\partial ^2 u \over \partial y ^2}\label{e40}\\
0&=&{\partial u \over \partial x}+{\partial v \over \partial y}\label{e41}
\end{eqnarray}

Together with boundary conditions $0=u=v$ at $y=0$ and $u=U(x,t)$ at $y=\infty$.\\

For 

\begin{eqnarray}\label{e42}
u=\sum_{i=0}^{\infty}\sum_{j=1}^{\infty}\sum_{k=0}^{\infty}A_{i,j,k}x^iy^jt^k\  
\end{eqnarray}

and 

\begin{eqnarray}\label{e44}
v=\sum_{i=0}^{\infty}\sum_{j=1}^{\infty}\sum_{k=0}^{\infty}B_{i,j,k}x^iy^jt^k\ 
\end{eqnarray}\\

a simplified equivalent version of the general defining equations for the coefficients of the solution of the general form of the Navier Stokes equations arises. It being noted that for $u$ and $v$, the starting value for the index $j$ is \emph{in fact $1$ and not $0$}, as a result of the fact that $u=0$ at $y=0$ as given above. (For $y=0$ we obtain $0=B_{i,0,k}=A_{i,0,k}$ as all other terms vanish in $u$ and $v$.)\\

Further, as noted by Schlichting as pressure is impressed\cite[page~110]{r32} upon the boundary layer by the external flow, the pressure gradient is \emph{fully defined} and given by the following expression relating to distal flow $U(x,t)$ being

\begin{equation}
{\partial U \over \partial t}+U{\partial U \over \partial x}=-{1 \over \rho}{\partial P \over \partial x}\label{e45}
\end{equation}

After substitution of $u$ and $v$ as defined by equations (\ref{e42}) and (\ref{e44}) into Prandlt's boundary layer equations (\ref{e40}) and \eqref{e41} and then using \eqref{e45} and applying the same methods as used to resolve the nonlinear terms in PVI and in the general form of the Navier Stokes equation, the following defining equation is obtained

\begin{equation}\label{e46}
\begin{array}{l}
\sum_{i=0}^{\infty}\sum_{j=1}^{\infty}\sum_{k=0}^{\infty}A_{i,j,k+1}(k+1)x^iy^jt^k\\
+\sum_{i=0}^{\infty}\sum_{j=2}^{\infty}\sum_{k=0}^{\infty}\sum_{p=0}^{i}\sum_{q=1}^{j-1}\sum_{r=0}^{k}A_{p,q,r}A_{i-p+1,j-q,k-r}(i-p+1)x^iy^jt^k\\
+\sum_{i=0}^{\infty}\sum_{j=1}^{\infty}\sum_{k=0}^{\infty}\sum_{p=0}^{i}\sum_{q=1}^{j}\sum_{r=0}^{k}B_{p,q,r}A_{i-p,j-q+1,k-r}(j-q+1)x^iy^jt^k\\
=\\
\sum_{i=0}^{\infty}\sum_{k=0}^{\infty}U_{i,k+1}(k+1)x^it^k\\
+\sum_{i=0}^{\infty}\sum_{k=0}^{\infty}\sum_{p=0}^{i}\sum_{q=0}^{k}U_{p,q}U_{i-p+1,k-q}(i-p+1)x^it^k\\
\nu\sum_{i=0}^{\infty}\sum_{j=0}^{\infty}\sum_{k=0}^{\infty}A_{i,j+2,k}(j+2)(j+1)x^iy^jt^k
\end{array}
\end{equation}\\

where $U(x,t)=\sum_{i=0}^{\infty}\sum_{k=0}^{\infty}U_{i,k}x^{i}t^{k}$.\\

Equation (\ref{e46}) can be further resolved by noting that only under the condition $j\ge 2$ does a general defining equation for the coefficients result. Accordingly two special subsidiary equations result for the special cases of $j=0$ and $j=1$ being 

\begin{eqnarray}\label{e48}
0=U_{i,k+1}(k+1)+\sum_{p=0}^{i}\sum_{q=0}^{k}U_{p,q}U_{i-p+1,k-q}(i-p+1)+2\nu A_{i,2,k}
\end{eqnarray}

for $j=0$ and $i,k\ge 0$\\

and\\

\begin{eqnarray*}
A_{i,1,k+1}(k+1)+\sum_{p=0}^{i}\sum_{r=0}^{k}B_{p,1,r}A_{i-p,1,k-r}=\nu A_{i,3,k}(3)(2)
\end{eqnarray*}

which can also be written as \\

\begin{eqnarray}\label{e50}
A_{i,3,k}={1 \over 6\nu}[A_{i,1,k+1}(k+1)+\sum_{p=0}^{i}\sum_{r=0}^{k}B_{p,1,r}A_{i-p,1,k-r}]
\end{eqnarray}

for $j=1$ and $i,k \ge 0$.\\

Equations (\ref{e48}) and (\ref{e50}) define $A_{i,2,k}$ and $A_{i,3,k}$
respectively. Similarly, the defining equation for $A_{i,j+2,k}$ for $j \ge 2$ and $i,k \ge 0$ can be extracted, by simple equation of coefficients, from (\ref{e46}) to yield 

\begin{eqnarray}\label{e52}
A_{i,j+2,k}&=&{1\over \nu (j+1)(j+2)}\Big(A_{i,j,k+1}(k+1)+\sum_{p=0}^{i}\sum_{q=1}^{j-1}\sum_{r=0}^{k}(i-p+1)A_{p,q,r}A_{i-p+1,j-q,k-r}\nonumber\\
&&+\sum_{p=0}^{i}\sum_{q=1}^{j}\sum_{r=0}^{k}B_{p,q,r}A_{i-p,j-q+1,k-r}(j-q+1)\Big)
\end{eqnarray}

for $j \ge 2$, $i,k \ge 0$.\\

Equations (\ref{e48}), (\ref{e50}) and (\ref{e52}) can be further simplified using relations derived from the continuity equation as follows.\\

By substituting $u=\sum_{i=0}^{\infty}\sum_{j=1}^{\infty}\sum_{k=0}^{\infty}A_{i,j,k}x^iy^jt^k\ $ and $v=\sum_{i=0}^{\infty}\sum_{j=1}^{\infty}\sum_{k=0}^{\infty}B_{i,j,k}x^iy^jt^k\ $ into the continuity equation \eqref{e41} the following equation is obtained

\begin{eqnarray}\label{e54}
0=\sum_{i=0}^{\infty}\sum_{j=1}^{\infty}\sum_{k=0}^{\infty}A_{i+1,j,k}(i+1)x^iy^jt^k+\sum_{i=0}^{\infty}\sum_{j=0}^{\infty}\sum_{k=0}^{\infty}B_{i,j+1,k}(j+1)x^iy^jt^k.  
\end{eqnarray}\\

Upon setting $j=0$ in (\ref{e54}) we  obtain-

\begin{eqnarray}\label{e56}
0=B_{i,1,k}
\end{eqnarray}

Further, the general defining equation from equation (\ref{e54}) for $i,k \ge 0$, $j \ge 1$ is 

\[0=A_{i+1,j,k}(i+1)+B_{i,j+1,k}(j+1)\] which can be written in a slightly more useful form, for computational purposes, as

\begin{eqnarray}\label{e58}
B_{i,j+1,k}=-{i+1 \over j+1}A_{i+1,j,k}
\end{eqnarray}

From equation (\ref{e56}) the defining equation for $A_{i,3,k}$, (\ref{e50}) becomes simply-

\begin{eqnarray}\label{e60}
A_{i,3,k}={1 \over 6 \nu}[A_{i,1,k+1}(k+1)]
\end{eqnarray}

for $i,k \ge 0$.

From equation (\ref{e58}) we obtain

\[ B_{p,q,r}=-{p+1 \over q}A_{p+1,q-1,r}\]

which when substituted into equation (\ref{e52}) yields

\begin{eqnarray}
A_{i,j+2,k}&=&{1\over \nu (j+1)(j+2)}\Big(A_{i,j,k+1}(k+1)\nonumber\\
&&+\sum_{p=0}^{i}\sum_{q=1}^{j-1}\sum_{r=0}^{k}(i-p+1)A_{p,q,r}A_{i-p+1,j-q,k-r}\label{e62}\\
&&-\sum_{p=0}^{i}\sum_{q=2}^{j}\sum_{r=0}^{k}{(p+1)(j-q+1)\over q}A_{p+1,q-1,r}A_{i-p,j-q+1,k-r}\Big)\nonumber
\end{eqnarray}

for $i,k \ge 0$ and $j \ge 2$.\\

Not only are the coefficients in the expansion of $u$ and $v$ now fully defined but the coefficients {\em can also be systematically constructed} from the above results. \\

The above results can be summarized as follows\\

$A_{i,2,k}$ for $i,k \ge 0$ is fully defined by equation (\ref{e48}).\\

$A_{i,3,k}$ for $i,k \ge 0$ is fully defined by equation (\ref{e60}).\\

$A_{i,j+2,k}$ for $i,k \ge 0$ $j \ge 2$ is fully defined by equation (\ref{e62}).\\

From the above it is seen that the only remaining unknown is $A_{i,1,k}$. However,  $A_{i,1,k}$ is fully and completely defined by way of matching $u$ against the external flow. In principle the match against the external flow occurs in the limit as $y \to \infty$, of course in practise the match occurs over a dimension, that is small in comparison with the reference length $L$. Accordingly, the solution $u$ can be matched against the external flow at a finite distance from the body surface so as to define $A_{i,1,k}$.\\

The above process will fully define $u$. It is also noted that $v$ is also fully defined since $B_{i,j+1,k}=-((i+1)/(j+1))A_{i+1,j,k}$  for $i,k \ge 0$, $j \ge 1$ from equation (\ref{e58}). Accordingly, from equation (\ref{e58}) we can write 

\begin{eqnarray}\label{e64}
B_{i,j,k}=-{i+1 \over j}A_{i+1,j-1,k}
\end{eqnarray}

for $j \ge 2$ which together with equations (\ref{e48})-(\ref{e62}) fully define $B_{i,j,k}$ (where $0=B_{i,0,k}=B_{i,1,k}$ for $i,k \ge 0$).\\

A few further simple adjustments will be made so as to put $u$ and $v$ in a form that is more suitable for calculation.\\

Equation (\ref{e62}) can be written, by letting $j=j'-2$ and subsequently dropping primes as 

\begin{eqnarray}\label{e66}
A_{i,j,k}&=&{1\over \nu (j-1)(j)}\Big(A_{i,j-2,k+1}(k+1)+\sum_{p=0}^{i}\sum_{q=1}^{j-3}\sum_{r=0}^{k}(i-p+1)A_{p,q,r}A_{i-p+1,j-2-q,k-r}\nonumber\\
&&-\sum_{p=0}^{i}\sum_{q=2}^{j-2}\sum_{r=0}^{k}{(p+1)(j-1-q)\over q}A_{p+1,q-1,r}A_{i-p,j-1-q,k-r}\Big)
\end{eqnarray}

for $i,k \ge 0$ and $j 
\ge 4$.\\

Further, so as to use equations (\ref{e48}), (\ref{e60}) and (\ref{e66}) to best effect for purposes of computation, $u$ will be written as

\begin{eqnarray}
u&=&\sum_{i=0}^{\infty}\sum_{j=1}^{\infty}\sum_{k=0}^{\infty}A_{i,j,k}x^{i}y^{j}t^{k}\nonumber\\
&=&\sum_{i=0}^{\infty}\sum_{k=0}^{\infty}A_{i,1,k}x^{i}yt^{k}+\sum_{i=0}^{\infty}\sum_{k=0}^{\infty}A_{i,2,k}x^{i}y^{2}t^{k}+\sum_{i=0}^{\infty}\sum_{k=0}^{\infty}A_{i,3,k}x^{i}y^{3}t^{k}\nonumber\\
&&+\sum_{i=0}^{\infty}\sum_{j=4}^{\infty}\sum_{k=0}^{\infty}A_{i,j,k}x^{i}y^{j}t^{k}\nonumber\\
&=&\sum_{i=0}^{\infty}\sum_{k=0}^{\infty}\Big[A_{i,1,k}y+A_{i,2,k}y^2+A_{i,3,k}y^3\Big]x^it^k
+\sum_{i=0}^{\infty}\sum_{j=4}^{\infty}\sum_{k=0}^{\infty}A_{i,j,k}x^{i}y^{j}t^{k}\nonumber
\end{eqnarray}

Equation (\ref{e66}) and a similar equation for $v$ derived by using equation (\ref{e64}) can be used to construct the solutions to Prandtl's boundary layer equations. These methods have been used to generate solutions which have then been substituted back into the boundary layer equations by using Waterloo Maple Software to verify the validity of the expansions. It has then been confirmed, that the solutions constructed to increasing orders of magnitude, lead to remainders (residuals) {\em that systematically vanish} as the order of expansion is increased.\\

\section{Conclusions}

It is seen that the above methods can be used to construct exact analytic solutions to many differential equations (ordinary and partial) that are likely to be encountered in practise. Nevertheless, many challenges remain for further research. The use of exact analytic expansions in numerical analytic continuation as outlined in Davis. Additional challenges of a non-trivial nature also include the desirability of constructing special solutions to differential equations in terms of classical transcendents (as noted above in relation to the work conducted by mathematicians in Belorussia (Lukashevich and Yablonski) and also by other contributors including Fokas and Ablowitz on Painlev\'{e} transcendents). \\

In closing, it has to be acknowledged that it is completely impossible in a paper of this size to do justice or pay tribute to the many important works worthy of further consultation in the areas of fluid dynamics and special functions. However, two works, by two leading authorities in their fields, provide a starting point, at the very least, for purposes of further investigation. In the area of fluid dynamics an excellent classical work that provides a thorough grounding in the principles of fluid dynamics is G.\ K.\ Batchelor's {\em An Introduction to Fluid Dynamics}\cite{r71} and in the area of the Painlev\'{e} transcendents, P.\ A.\ Clarkson's review of the Painlev\'{e} transcendents in the {\em NIST Handbook of Mathematical Functions}\cite{r27} also provides an excellent starting point for further research.\\

\section*{Acknowledgments}

This paper is dedicated to the memory of the late Professor P.\ L.\ Sachdev, formerly Professor of Mathematics, Indian Institute of Science, Bangalore, India.  \\ 

The author also gratefully acknowledges the support of other persons who have made his work on series solutions possible over the past 20 years including, without limitation,\\ Prof.\ T.\ Marchant of the University of Wollongong's School of Mathematics and Applied Statistics; Prof.\ A.\ Shannon and Dr G.\ Smith of the University of Technology Sydney's School of Mathematics for their support and encouragement during the early 1990s.\\

\end{document}